\documentclass[12pt]{amsart}

\newtheorem{theorem}{Theorem}[section]
\newtheorem{lemma}[theorem]{Lemma}

\newtheorem{corollary}[theorem]{Corollary}

\theoremstyle{definition}

\theoremstyle{remark}

\numberwithin{equation}{section}

\def\ot{\otimes}

\def\pa{\partial}
\def\pab{\bar{\partial}}

\def\om{\omega}

\def\we{\wedge}

\def\rhob{\bar{\rho}}

\def\rhob{\bar{\rho}}

\def\imi{\sqrt{-1}}

\newcommand\Lam{\Lambda}

\newcommand\bbc{\mathbb{C}}

\newcommand\ra{\longrightarrow}

\begin{document}
	
	\title[Dual Forms]{Dual Forms of the Squares of Nijenhuis Tensor}
	
	\author[Jun Ling]{Jun \underline{LING}}
	\address{Department of Mathematics, Utah Valley University, Orem, Utah 84058}
	\email{lingju@uvu.edu}
	
	
	
	\date{}
	
	
	
	\begin{abstract}
		We express the  dual forms of squares of Nijenhuis tensor  in terms of the second order component derivatives
		of the exterior derivative on differential forms and  give new vanishing results for the squares of Nijenhuis tensor. 
	\end{abstract}
	
	\maketitle
	
	\section{Introduction}\label{sec-intro}
	A complex manifold is  a differentiable manifold,  for the holomorphic structure on the complex manifold gives
	a differentiable structure. The converse is not true in general. It is known n-spheres ($n\not=2,6$) are not complex manifold (Borel and Serre \cite{B-S}), though they are differentiable manifolds. On the other hand, it is well-known that 2-sphere  is a differentiable manifold and a complex manifold as well.  One would like to know whether a given differential manifold is also a complex manifold so that the holomorphic structure induces the same differential structure as the one on the underlying differentiable manifold. 
	
	\vspace{0.2in}
	
	The existence of complex manifold structures on high dimensional (more than two) spheres has been studied extensively, cf. Hopf \cite{Hopf},
	Borel and Serre \cite{B-S}, 
	Ehresmann \cite{Ehr}, Kirchhoff \cite{Kir},
	Eckmann and Fr\"{o}licher \cite{E-F},  Ehressmann and Libermann \cite{E-L}, LeBrun \cite{LeBrun}, and etc, left with an open case of 6-sphere, cf.  Hirzebruch \cite{Hir} in 1954, Libermann \cite{Lib} in 1955 and
	Yau \cite{Yau} in 1990.
	
	\vspace{0.2in}
	
	A complex manifold has a complex structure that is an almost-complex structure on the underlying differentiable manifold. An almost-complex structure $J$ on a differentiable manifold is an
	endomorphism on the tangent bundle of the manifold with $J^2=-1$. Conversely
	an integrable almost-complex structure  on a differentiable manifold  is a complex structure that makes the differentiable manifold a complex manifold.  So for differentiable manifolds with almost-complex structures, say for open case 6-sphere, the problem becomes whether
	the existing almost-complex structures are integrable. 
	
	\vspace{0.2in}
	
	Newlander and Nirenberg \cite{NN} showed that an almost-complex structure  on a differentiable manifold is integrable if and only if  Nijenhuis tensor of the almost-complex structure  vanishes. Thus study of Nijenhuis tensor  and its vanishing conditions
	has been a focus, see for example the work of  Fr\"{o}licher and Nijenhuis
	\cite{F-N}, Newlander and Nirenberg \cite{NN},   Cirici and Wilson \cite{C-W}, and others (see also \cite{Ling}). The Nijenhuis tensor $N$ of an almost-complex structure $J$ on a differentiable manifold $M$ can be defined as follows.
	
	\begin{equation}\label{N-def}
	\left\{
	\begin{array}{l}
	N: \Gamma(TM)\times \Gamma(TM)\ra \Gamma(TM), \quad \forall X,Y\in \Gamma(TM),\\
	N(X,Y):=[JX,JY]-J[X,JY]-J[JX,Y]-[X,Y],
	\end{array}
	\right.
	\end{equation}
	where $[\cdot,\cdot]$ is the Lie bracket.
	
	\vspace{0.2in}
	
	One approach to  Nijenhuis tensor $N$ and its vanishing conditions is to study the "squares" of $N$, and vanishing results of those squares. The motivation  is to relate the Nijenhuis tensor to cohomology of almost-complex manifold, cf. \cite{C-W} and others. 
	Our work is in this direction.
	
	\vspace{0.2in}
		
\cite{Ling} considered the squares derived from $JN^2$ and vanishing results, in terms 
of  matrices of almost-complex structure $J$ and Nijenhuis tensor $N$ under local coordinate of the manifolds, where the products are compositions of maps on appropriate arguments. One of main vanishing results there  has been established  by  algebraic calculations of the sixteen terms of five factor products of the matrices of $J$ and its derivative matrices, thanks to  an extra factor $J$ in $JN^2$ and the identity $J^2=-1$.
However, that algebraic argument  would not produce necessary cancellations in establishing vanishing results 
for the squares derived from $N^2$, without an extra factor $J$.
Some new method needs to be developed to study vanishing results in this new case. In this paper we give the dual forms of the squares derived from $N^2$  first, then prove our new vanishing results by analytic method. We are able to give dual forms of the squares derived from the $JN^2$  in \cite{Ling} as well.
It is interesting that the dual forms for both the squares derived from $N^2$ and the ones derived from $JN^2$  are "conjugate" each other with their counterparts;  each square and its counterpart together produce a single term consisting of second order  component derivatives of the exterior derivative and  real objects through  addition and subtraction, similar to the situation that classic Euler formulas $e^{i\theta}=\cos(\theta)+i\sin(\theta)$ and $e^{-i\theta}=\cos(\theta)-i\sin(\theta)$ derive real valued cosine and sine functions $\cos\theta)$ and $\sin(\theta)$, see Theorems \ref{NN-N-thm} and \ref{NN-thm}. The dual forms of  the squares is an asset to the squares  since they are favored in relations with theories and machines in  Topology, and in Complex and Algebraic Geometry.

	\vspace{0.2in}
	
	This paper is organized as follows. We 
	state general settings and lemmas that establish relations among Nijenhuis tensor  and component derivatives of the exterior derivative in Section \ref{sec-preparation}. The dual forms are developed in two steps in two separate sections. Section \ref{sec-dual-form-1}  is the first step, where we give the first type of dual forms that express the squares of Nijenhuis tensor in terms of Nijenhuis tensor and Lie bracket. We then give new vanishing results.
	Section \ref{sec-dual-form-2} is the second step, where we present the second type of dual forms of  the squares  and vanishing results  in terms of component derivatives  of the exterior derivative.

	\vspace{0.2in}
	
	\section{Preparation}\label{sec-preparation}
	In this paper, we let $M$ be a $n=2m$-dimensional almost-complex manifold with almost-complex structure $J$. Note $J$ is an endomorphism on the tangent bundle $TM$ with $J^2=-1$. 
	Let $N$ be the Nijenhuis tensor on  $TM$, that is given by (\ref{N-def}).

	\vspace{0.1in}
	
	Define $N^2$ by as follows.
	\begin{equation}\label{N-sq-def}
	\left\{
	\begin{array}{l}
	N^2: \Gamma(TM)\times \Gamma(TM) \times \Gamma(TM) \ra \Gamma(TM),\\
	N^2(X,Z; Y)=N\{N(X,Z),Y\}, \ \forall X,Y\in \Gamma(TM).
	\end{array}
	\right.
	\end{equation}
	
	We call $N^2$ and $JN^2$  the \textbf{strong squares}  of Nijenhuis tensor $N$.
	
	\vspace{0.1in}
	
	Since the eigenvalues $\pm\imi$ of $J_x$ for each $x\in M$ are pure imaginary numbers, we need to complexify objects.
	
	For every $x\in M$, let $T_xM^{\bbc}:=T_xM\ot \bbc$ be the complexification of the tangent space $T_xM$ of $M$ at $x$.
	It is known that
	\[
	T_xM^{\bbc}=T_x^{1,0}M\oplus T_x^{0,1}M,
	\]
	where $T_x^{1,0}M$ and $T_x^{0,1}M$ are eigenspaces of eigenvalue $\imi$ and eigenvalue $-\imi$ of $J_x$,  respectively.
	
	\vspace{0.1in}
	
	Denote the complexified cotangent bundle by  $T^*M^{\bbc}=T^*M\ot \bbc$,\\
	
	\vspace{0.1in}
	
	Taking the duals and wedge products, the above decomposition of complexified tangent bundle  and
	complexified cotangent bundle of $M$ gives the decomposition
	\[
	\Lambda^k=\oplus_{p+q=k}\Lambda^{p,q}
	\]
	where $\Lambda^{p,q}$ consists of forms with bigrade $(p,q)$. At $x\in M$,
	\[
	\Lambda_x^{p,q}=\wedge_p(T_x^*M^\bbc)^{1,0}\ot \wedge_q(T_x^*M^\bbc)^{0,1}.
	\]

	\vspace{0.1in}

	$\Lambda^{0,0}$, $\Lambda^{1,0}$, and $\Lambda^{0,1}$ generate $\Lambda^{p,q}$ for all non-negative integers $p,q$.
	Let $d$ be the exterior derivative on differentiable forms.
	It is known that $d^2=0$ and 
	\[
	d(\Lambda^{p,q})\subseteq \Lambda^{p-1,q+2}\oplus \Lambda^{p,q+1}\oplus \Lambda^{p+1,q}\oplus \Lambda^{p+2,q-1}.
	\]
	For a form $\om$ we use $\pi^{p,q}\om$ for the projection of $\om$ on $\Lambda^{p,q}$ and write 
	\[
	d=\rhob+\pab+\pa+\rho,
	\]
	where \\
	$\rhob=d^{-1,2}$ is the bidegrees $(-1,2)$ component ,\\
	$\rho=d^{2,-1}$ the bidegrees $(2,-1)$ component,\\ 
	$\pab=d^{0,1}$ the bidegrees $(0,1)$ component,\\
	$\pa=d^{1,0}$ the bidegrees $(1,0)$ component of $d$, respectively.

	Note that Leibnize rule applies to exterior derivative and each component derivative. In particular
	\[
	\rhob(\om\we \theta)=\rhob\om\we\theta+(-1)^{\deg(\om)}\om\we\rhob\theta,
	\]
	\[
	\rho(\om\we \theta)=\rho\om\we\theta+(-1)^{\deg(\om)}\om\we\rho\theta.
	\]

	For function $f$ on $M$, $\rho f=0$, $\rhob f=0$ and above rules imply that  $\rho$ and $\rhob$ are linear over functions, so  maps $\rho$ and $\rhob$ may restrict to fibers: for $x\in M$,
	\[
	\rho_x:\Lam_x^{p,q}\ra \Lam_x^{p+2,q-1},\quad \rhob_x:\Lam_x^{p,q}\ra \Lam_x^{p-1,q+2}.
	\]
	We will omit the subscript $x$ often.
	\vspace{0.1in}

	For a vector $X$ we use $\pi_{0,1}X$, $\pi_{1,0}X$ for the (0,1)-component and (1,0)-comonent of $X$, respectively.
	
	\vspace{0.1in}
	
	We present some Lemmas needed for the formulations in the next section. Most of them could be well-known.
	
	\vspace{0.1in}
	
	The first lemma describe the  Nijenhuis tensor $N$  on complexified tangent bundle of the manifold
	$TM^{\bbc}=T_{1,0}M\oplus T_{0,1}M$.
	
	\begin{lemma}\label{Lie-N} Let $X,Y$ be real vector fields, then we have
		\begin{equation}\label{10bracket-01N}
		[\pi_{1,0}X,\pi_{1,0}Y]=-\frac14\pi_{0,1}N(X,Y),
		\end{equation}
		\begin{equation}\label{01bracket-10N}
		[\pi_{0,1}X,\pi_{0,1}Y]=-\frac14\pi_{1,0}N(X,Y).
		\end{equation}
	\end{lemma}
	\proof
	\[
	[\pi_{1,0}X,\pi_{1,0}Y]
	\]
	\[
	=\frac14[X-\imi JX,Y-\imi JY]
	\]
	\[
	=\frac14\{[X,Y]-\imi[X,JY]-\imi[JX,Y]-[JX,JY] \}
	\]
	\[
	=-\frac14\{[JX,JY]+\imi[X,JY]+\imi[JX,Y]-[X,Y] \}.
	\]
	
	So
	
	\[
	\pi_{0,1}[\pi_{1,0}X,\pi_{1,0}Y]
	\]
	\[
	=-\frac18\{[JX,JY]+\imi[X,JY]+\imi[JX,Y]-[X,Y] 
	\]
	\[
	+\imi J[JX,JY]+\imi J\imi[X,JY]+\imi J\imi[JX,Y]-\imi J[X,Y] \}
	\]
	\[
	=-\frac18\{[JX,JY]+\imi[X,JY]+\imi[JX,Y]-[X,Y] 
	\]
	\[
	+\imi J[JX,JY]-J[X,JY]-J[JX,Y]-\imi J[X,Y] \}
	\]
	\[
	=-\frac18\{[JX,JY]-J[X,JY]-J[JX,Y]-[X,Y] 
	\]
	\[
	+\imi J[JX,JY]+\imi[X,JY]+\imi[JX,Y]-\imi J[X,Y] \}
	\]
	\[
	=-\frac18\{N(X,Y)
	+\imi JN(X,Y) \}
	\]
	\[
	=-\frac14\{\pi_{0,1}N(X,Y)\}.
	\]
	Therefore
	\[
	[\pi_{1,0}X,\pi_{1,0}Y]=-\frac14\pi_{0,1}N(X,Y).
	\]
	
	Similarly,
	
	\[
	[\pi_{0,1}X,\pi_{0,1}Y]
	\]
	\[
	=\frac14[X+\imi JX,Y+\imi JY]
	\]
	\[
	=\frac14\{[X,Y]+\imi[X,JY]+\imi[JX,Y]-[JX,JY] \}
	\]
	\[
	=-\frac14\{[JX,JY]-\imi[X,JY]-\imi[JX,Y]-[X,Y] \}.
	\]
	
	So
	
	\[
	\pi_{1,0}[\pi_{0,1}X,\pi_{0,1}Y]
	\]
	\[
	=-\frac18\{[JX,JY]-\imi[X,JY]-\imi[JX,Y]-[X,Y] 
	\]
	\[
	-\imi J[JX,JY]+\imi J\imi[X,JY]+\imi J\imi[JX,Y]+\imi J[X,Y] \}
	\]
	\[
	=-\frac18\{[JX,JY]-\imi[X,JY]-\imi[JX,Y]-[X,Y] 
	\]
	\[
	-\imi J[JX,JY]-J[X,JY]-J[JX,Y]+\imi J[X,Y] \}
	\]
	\[
	=-\frac18\{[JX,JY]-J[X,JY]-J[JX,Y]-[X,Y] 
	\]
	\[
	-\imi J[JX,JY]-\imi[X,JY]-\imi[JX,Y]+\imi J[X,Y] \}
	\]
	\[
	=-\frac18\{[JX,JY]-J[X,JY]-J[JX,Y]-[X,Y] 
	\]
	\[
	-\imi J([JX,JY]-J[X,JY]-J[JX,Y]-[X,Y]) \}
	\]
	\[
	=-\frac18\{N(X,Y)
	-\imi JN(X,Y) \}
	\]
	\[
	=-\frac14\{\pi_{1,0}N(X,Y)\}.
	\]
	
	Therefore
	
	\[
	[\pi_{0,1}X,\pi_{0,1}Y]=-\frac14\pi_{1,0}N(X,Y).
	\]

	\qed
	
	\vspace{0.1in}
	
	The following lemma describe the relation of Nijenhuis tensor $N$ with components $\rho,\rhob,\pa, \pab$ of exterior derivative $d$.
	
	\begin{lemma}\label{N-mu-lem}
		Let $f$ be a smooth function on $M$, and $X$ and $Y$ be real tangent vectors. Then we have 
		the following equations\\
		(a)
		\[
		\pab f\big(N(X,Y)\big)=-4(\pa^2 f)(X,Y).
		\]
		(b)
		\[
		\pa f\big(N(X,Y)\big)=-4(\pab^2 f)(X \we Y).
		\]
		(c) If	$\om$ is a (0,1)-form then we have
		\[
		\om(N(X,Y))=4(\rho\om)(X,Y).
		\]
		(d) If	$\theta$ is a (1,0)-form then we have
		\[
		\theta(N(X,Y))=4(\rhob\theta)(X,Y).
		\]
	\end{lemma}
	Cirici and Wilson \cite{C-W} gave results (c) and (d).
	
	For convenience, we present a proof for all (a)-(d).
	
	\proof
	For real vector fields $X$ and $Y$, Let $\pi_{1,0}X=X_{1,0}=\frac12(X-\imi JX)$ and $\pi_{0,1}X=X_{0,1}=\frac12(X+\imi JX)$. Then
	\[
	\pab^2f(X,Y)
	\]
	\[
	=\{\pi^{0,1}d(\pab f)\}\big((X_{1,0}+X_{0,1})\we(Y_{1,0}+Y_{0,1})\big)
	\]
	\[
	=\{\pi^{0,1}d(\pab f)\}(X_{1,0}\we Y_{1,0})
	\]
	\[
	+\{\pi^{0,1}d(\pab f)\}(X_{1,0}\we Y_{0,1})
	\]\[
	+\{\pi^{0,1}d(\pab f)\}(X_{0,1}\we Y_{1,0})
	\]
	\[
	+\{\pi^{0,1}d(\pab f)\}(X_{0,1}\we Y_{0,1})
	\]
	\[
	=\{\pi^{0,1}d(\pab f)\}(X_{0,1}\we Y_{0,1})
	\]
	\[
	=\{\pi^{0,1}d(\pab f)\}(\pi_{0,1}X\we \pi_{0,1}Y)
	\]
	\[
	=d(\pab f)(\pi_{0,1}X,\pi_{0,1}Y)
	\]
	\[
	=\pi_{0,1}X\big(\pab f(\pi_{0,1}Y)\big)-\pi_{0,1}Y\big(\pab f(\pi_{0,1}X)\big)
	-\pab f([\pi_{0,1}X,\pi_{0,1}Y]).
	\]
	
	By Lemma \ref{Lie-N},
	the above equation becomes
	\[\pi_{0,1}X\big((\pi_{0,1}Y)(f)\big)-\pi_{0,1}Y\big((\pi_{0,1}X)(f)\big)-\pab f(-\frac14\pi_{1,0}N(X,Y)).
	\]
	\[
	=\pi_{0,1}X\big((\pi_{0,1}Y)(f)\big)-\pi_{0,1}Y\big((\pi_{0,1}X)(f)\big)
	\]
	\[
	=[\pi_{0,1}X, \pi_{0,1}Y](f)
	\]
	\[
	=-\frac14\pi_{1,0}N(X,Y)(f)
	\]
	\[
	=-\frac14\pa f\big(N(X,Y)\big),
	\]
	where in the second equation to the last one, we used Lemma \ref{Lie-N} again.
	Therefore we have
	\[
	\pa f\big(N(X,Y)\big)=-4(\pab^2 f)(X,Y).
	\]
	
	\vspace{0.1in}
	
	Similarly,
	\[
	\pa^2f(X,Y)
	\]
	\[
	=\{\pi^{1,0}d(\pa f)\}\big((X_{1,0}+X_{0,1})\we(Y_{1,0}+Y_{0,1})\big)
	\]
	\[
	=\{\pi^{1,0}d(\pa f)\}(X_{1,0}\we Y_{1,0})
	\]
	\[
	=\{\pi^{1,0}d(\pa f)\}(\pi_{1,0}X\we\pi_{1,0}Y)
	\]
	\[
	=d(\pa f)(\pi_{1,0}X\we\pi_{1,0}Y)
	\]
	\[
	=\pi_{1,0}X\big(\pa f(\pi_{1,0}Y)\big)-\pi_{1,0}Y\big(\pa f(\pi_{1,0}X)\big)
	-\pa f([\pi_{1,0}X,\pi_{1,0}Y]).
	\]
	
	By Lemma \ref{Lie-N},
	the above equation becomes
	\[\pi_{1,0}X\big((\pi_{1,0}Y)(f)\big)-\pi_{1,0}Y\big((\pi_{1,0}X)(f)\big)-\pa f(-\frac14\pi_{0,1}N(X,Y)).
	\]
	\[
	=\pi_{1,0}X\big((\pi_{1,0}Y)(f)\big)-\pi_{1,0}Y\big((\pi_{1,0}X)(f)\big)
	\]
	\[
	=[\pi_{1,0}X, \pi_{1,0}Y](f)
	\]
	\[
	=-\frac14\pi_{0,1}N(X,Y)(f)
	\]
	\[
	=-\frac14\pab f\big(N(X,Y)\big),
	\]
	where in the second equation to the last one, we used Lemma \ref{Lie-N} again.
	Therefore we have
	\[
	\pab f\big(N(X,Y)\big)=-4(\pa^2 f)(X,Y).
	\]

	\vspace{0.1in}
	
	For  $(0,1)$-form $\om$ and real vector fields
	$X$ and $Y$, we compute
	\[
	(\rho\om)(X,Y)=\{\pi^{2,0}(d\om)\}(X,Y)
	\]
	\[
	=\{\pi^{2,0}(d\om)\}\big((X^{1,0}+X^{0,1})\we (Y^{1,0}+Y^{0,1})\big)
	\]
	\[
	=\{\pi^{2,0}(d\om)\}(X_{1,0}\we Y_{1,0})
	\]
	\[
	+\{\pi^{2,0}(d\om)\}(X_{1,0}\we Y_{0,1})
	\]\[
	+\{\pi^{2,0}(d\om)\}(X_{0,1}\we Y_{1,0})
	\]\[
	+\{\pi^{2,0}(d\om)\}(X_{0,1}\we Y_{0,1})
	\]
	\[
	=\{\pi^{2,0}(d\om)\}(X_{1,0}\we Y_{1,0})
	\]
	\[
	=\{\pi^{2,0}(d\om)\}(\pi_{1,0}X\we \pi_{1,0}Y)
	\]
	\[
	=(d\om)(\pi_{1,0}X\we \pi_{1,0}Y)
	\]
	\[
	=\pi_{1,0}X\big(\om(\pi_{1,0}Y)\big)-\pi_{1,0}Y\big(\om(\pi_{1,0}X)\big)-\om([\pi_{1,0}X,\pi_{1,0}Y])
	\]
	\[
	=-\om([\pi_{1,0}X,\pi_{1,0}Y])
	\]
	\[
	=-\om(\pi_{0,1}[\pi_{1,0}X,\pi_{1,0}Y]).
	\]
	
	By  Lemma \ref{Lie-N}:
	\[
	[\pi_{1,0}X,\pi_{1,0}Y]=-\frac14\pi_{0,1}N(X,Y).
	\]
	
	Therefore,
	\[
	(\rho\om)(X,Y)=\frac14\om(\pi_{0,1}N(X,Y)),
	\]
	\[
	=\frac14\om(N(X,Y)).
	\]
	
	Similarly, for $(1,0)$-form $\theta$ and real vector fields $X$ and $Y$, we have the following.
	\[
	(\rhob\theta)(X,Y)=\{\pi^{0,2}(d\theta)\}(X,Y)
	\]
	\[
	=\{\pi^{0,2}(d\theta)\}\big((X_{1,0}+X_{0,1})\we(Y_{1,0}+Y_{0,1})\big)
	\]
	\[
	=\{\pi^{0,2}(d\theta)\}(X_{0,1},Y_{0,1})
	\]
	\[
	=\{\pi^{0,2}(d\theta)\}(\pi_{0,1}X,\pi_{0,1}Y)
	\]
	\[
	=(d\theta)(\pi_{0,1}X,\pi_{0,1}Y)
	\]
	\[
	=\pi_{0,1}X\big(\theta(\pi_{0,1}Y)\big)-\pi_{0,1}Y\big(\theta(\pi_{0,1}X)\big)
	-\theta([\pi_{0,1}X,\pi_{0,1}Y])
	\]
	\[
	=-\theta([\pi_{0,1}X,\pi_{0,1}Y])
	\]
	\[
	=-\theta(\pi_{1,0}[\pi_{0,1}X,\pi_{0,1}Y]).
	\]
	\[
	=\frac14\theta(\pi_{1,0}N(X,Y)),
	\]
	\[
	=\frac14\theta(N(X,Y)),
	\]
	where by Lemma \ref{Lie-N},
	\[
	[\pi_{0,1}X,\pi_{0,1}Y]=-\frac14\pi_{1,0}N(X,Y).
	\]
	
	Therefore for $(0,1)$-form $\theta$
	\[
	(\rhob\theta)(X,Y)=\frac14\theta(N(X,Y)).
	\]
	
	\qed

	\vspace{0.2in}

	\section{Type I Dual Forms of the Squares of Nijenhuis Tensor}\label{sec-dual-form-1}

	In this section we present a dual forms of the strong squares $N^2$  and $JN^2$ of Nijenhuis tensor in terms of Nijenhuis tensor itself  and Lie bracket, prove a vanishing result and identity.
	
	\begin{theorem}\label{NN-N-thm} We have the following results.
		
		For (0,1)-form $\om$,
		\begin{equation}\label{01-NN-N-eq}
		\om\big( N^2(X,Z; Y))\big)=\om\big([N(X, Z), Y]\big)
		\end{equation}
		and
		\[
		\om\big(J N^2(X,Z; Y)\big)=-\imi\om\big([N(X, Z), Y]\big).
		\]

		\vspace{0.1in}
		
		For (1,0)-form $\theta$,
		\begin{equation}\label{10-NN-N-eq}
		\theta\big( N^2(X,Z; Y)\big)=\theta\big([N(X, Z), Y]\big)
		\end{equation}
		and
		\[
		\theta\big( JN^2(X,Z; Y)\big)=\imi\theta\big([N(X, Z), Y]\big).
		\]
		
		\vspace{0.1in}
		
		Therefore one can extend $N^2$ and $JN^2$ to derivations on exterior algebra of differential forms.	\\
		
		We also have the following.
		
		\vspace{0.1in}
		
		For real form $\zeta$,
		\begin{equation}\label{real-NN-N-eq}
		\zeta\big(N^2(X,Z; Y)\big)=\zeta\big([N(X, Z), Y]\big)
		\end{equation}
		and
		\[
		\zeta\big(JN^2(X,Z; Y)\big)=(J\zeta)\big([N(X, Z), Y]\big).
		\]
		\vspace{0.1in}

	\end{theorem}
	
	\proof
	For $(0,1)$ form $\om$, real vectors $X,Z$ and $Y$,  we apply Lemma \ref{N-mu-lem}. Then
	\[
	\om\big(N^2(X,Z; Y)\big)
	\]
	\[
	=\om\Big(N\big(N(X, Z), Y\big)\Big)
	\]
	\[
	=4(\rho \om)\Big(N(X, Z), Y\Big)
	\]
	\[
	=4(\pi^{2,0}d \om)\Big(N(X, Z), Y\Big)
	\]
	\[
	=4(\pi^{2,0}d \om)\Big(\pi_{1,0}N(X, Z), \pi_{1,0}Y\Big)
	\]
	\[
	=4(d \om)\Big(\pi_{1,0}N(X, Z), \pi_{1,0}Y\Big)
	\]
	\[
	=4\{\pi_{1,0}N(X, Z)\om(\pi_{1,0}Y)-\pi_{1,0}Y\om(\pi_{1,0}N(X, Z))-\om([\pi_{1,0}N(X, Z), \pi_{1,0}Y])\}
	\]
	\[
	=-4\om([\pi_{1,0}N(X, Z), \pi_{1,0}Y])
	\]
	\[
	=\om(\pi_{0,1}[N(X, Z), Y])
	\]
	\[
	=\om([N(X, Z), Y])
	\]

	\[
	\om\big(N^2(X,Z; Y)\big)=\om\big([N(X, Z), Y]\big)
	\]

	
	\[
	\theta\big(N^2(X,Z; Y)\big)
	\]
	\[
	=\theta\Big(N\big(N(X, Z), Y\big)\Big)
	\]
	\[
	=4(\rhob \theta)\Big(N(X, Z), Y\Big)
	\]
	\[
	=4(\pi^{0,2}d \theta)\Big(N(X, Z), Y\Big)
	\]
	\[
	=4(\pi^{0,2}d \theta)\Big(\pi_{0,1}N(X, Z), \pi_{0,1}Y\Big)
	\]
	\[
	=4(d \theta)\Big(\pi_{0,1}N(X, Z), \pi_{0,1}Y\Big)
	\]
	\[
	=4\{\pi_{0,1}N(X, Z)\theta(\pi_{0,1}Y)-\pi_{0,1}Y\theta(\pi_{0,1}N(X, Z))-\theta([\pi_{0,1}N(X, Z), \pi_{0,1}Y])\}
	\]
	\[
	=-4\theta([\pi_{0,1}N(X, Z), \pi_{0,1}Y])
	\]
	\[
	=\theta(\pi_{1,0}[N(X, Z), Y])
	\]
	\[
	=\theta([N(X, Z), Y])
	\]
	
	\[
	\theta\big(N^2(X,Z; Y)\big)=\theta([N(X, Z), Y])
	\]
	
	By Lemma \ref{J-through-01-10-form-lem}
	\[
	\om\big(N^2(X,Z; Y)\big)=-\imi\om\Big(N\big(N(X, Z), Y\big)\Big)=-\imi\om\big([N(X, Z), Y]\big),
	\]
	\[
	\theta\big(N^2(X,Z; Y)\big)=\imi\theta\Big(N\big(N(X, Z), Y\big)\Big)=\imi\theta\big([N(X, Z), Y]\big).
	\]
	The proof for $\om$ and $\theta$ is done.
	
	For real form $\zeta=\frac12(\zeta+\imi J\zeta)+\frac12(\zeta-\imi J\zeta)$
	\[
	\frac12(\zeta+\imi J\zeta)\big(N^2(X,Z; Y)\big)=\frac12(\zeta+\imi J\zeta)\big([N(X, Z), Y]\big)
	\]
	\[
	\frac12(\zeta-\imi J\zeta)\big(N^2(X,Z; Y)\big)=\frac12(\zeta-\imi J\zeta)\big([N(X, Z), Y]\big).
	\]
	Add the above two equations and add the two equations respectively, 
	\[
	\frac12(\zeta+\imi J\zeta)\big(N^2(X,Z; Y)\big)=-\imi\frac12(\zeta+\imi J\zeta)\big([N(X, Z), Y]\big)
	\]
	\[
	\frac12(\zeta-\imi J\zeta)\big(N^2(X,Z; Y)\big)=\imi\frac12(\zeta-\imi J\zeta)\big([N(X, Z), Y]\big)
	\]
	we get the last two equations for real form $\zeta$.
	\qed

	For real vectors $X,Z,Y,W$ and their duals $X^*,Z^*, Y^*, W^*$, we define the following functionals.
	\[
	K(X,Z,Y,W):=
	\]
	\[
	\frac14W^*([N(X, Z), Y])
	\]
	\[+\frac14W^*([N(Y, Z), X])
	\]
	\[+\frac14Z^*([N(X, W), Y])
	\]
	\[+\frac14Z^*([N(Y, W), X]),
	\]
	
	\[
	L(X,Z,Y,W):=
	\]
	\[
	\frac14W^*(J[N(X, Z), Y])
	\]
	\[+\frac14W^*(J[N(Y, Z), X])
	\]
	\[+\frac14Z^*(J[N(X, W), Y])
	\]
	\[+\frac14Z^*(J[N(Y, W), X]),
	\]
	
	We now define the \textbf{intermediate squares} $\hbar$ and $\ell$ of Nijenhuis tensor $N$ by
	\[
	\hbar(X,Z):=Z^*([N(X, Z), X]),
	\]
	and
	\[
	\ell(X,Z):=Z^*([N(X, Z), X]).
	\]
	and 
	the \textbf{weak squares} $S$ and $T$ at $x\in M$ of Nijenhuis tensor $N$ by
	\[
	S=\sum_{i,k=1}^n\hbar(\pa_i,\pa_k)
	\]
	\[
	T=\sum_{i,k=1}^n\ell(\pa_i,\pa_k),
	\]
	
	where $\{x^i\}_{i=1}^{n=2m}$ be a local chart of $M$ at $x$, 
	and $\{dx^i\}_{i=1}^{n}$ be the dual of $\{\dfrac{\pa}{\pa x^i}\}_{i=1}^n:=\{\pa_i\}_{i=1}^n$
	Let
	\[S_i:=\sum_{k=1}^n\hbar(\pa_i,\pa_k),\quad\textup{for }1\le i\le n.\]
	
	\cite{Ling} showed vanishing result $T=0$. Now we give other two vanishing results.
	\begin{theorem}\label{S-vanishes-thm}
		\[
		S_i\equiv 0,\qquad \textup{for}\quad  i=1,2,\cdots, n
		\] 
		and
		\[
		S\equiv 0.
		\]
	\end{theorem}
	\proof
	Denote
	\[
	N_{ik}^r:=dx^r\big(N(\pa_i,\pa_k)\big).
	\]
	It is easy to know that
	\[
	N_{ik}^r=\sum_{p=1}^nJ_i^p\{\pa_pJ_k^r-\pa_kJ_p^r\}-J_k^p\{\pa_pJ_i^r-\pa_iJ_p^r\},
	\]
	\[
	N_{ik}^k=\sum_{p=1}^nJ_i^p\{\pa_pJ_k^k-\pa_kJ_p^k\}-J_k^p\{\pa_pJ_i^k-\pa_iJ_p^k\}
	\]
	for each $i$ and $k$ with $1\le i,k \le n$. So
	\[
	\sum_{k=1}^nN_{ik}^k=\sum_{p,k=1}^nJ_i^p\pa_pJ_k^k+\sum_{p,k=1}^n-J_i^p\pa_kJ_p^k
	+\sum_{p,k=1}^n-J_k^p\pa_pJ_i^k+\sum_{p,k=1}^nJ_k^p\pa_iJ_p^k
	\]
	\[
	:=I+ II+ III +IV
	\]

	The first term in the last equation is
	\[
	I=\sum_{p,k=1}^nJ_i^p\pa_pJ_k^k=\sum_{p,k=1}^n\pa_p\sum_kJ_k^k=0.
	\]

	The second term 
	\[
	II=\sum_{p,k=1}^n-J_i^p\pa_kJ_p^k
	\]
	\[
	=\sum_{p,k=1}^nJ_p^k\pa_kJ_i^p
	\]
	\[
	=\sum_{p,k=1}^nJ_k^p\pa_pJ_i^k
	\]
	\[
	=-III.
	\]
	which is the negative of the third term.
	The fourth term
	\[
	IV=\sum_{p,k=1}^nJ_k^p\pa_iJ_p^k
	\]
	\[
	=-\sum_p\sum_kJ_p^k\pa_iJ_k^p
	\]
	\[
	=-\sum_k\sum_pJ_k^p\pa_iJ_p^k
	\]
	\[
	=-IV.
	\]
	Therefore we have
	\begin{equation}\label{sum-N-ikk-is-0}
	\sum_kN_{ik}^k=0.
	\end{equation}

	On the other hand, by Theorem \ref{NN-N-thm}
	\[
	\hbar(\pa_i,\pa_k)=(dx^k)^*([N(\pa_i, \pa_k), \pa_i])
	\]
	\[
	=(dx^k)([N_{ik}^r\pa_r, \pa_i])
	\]
	\[
	=-(dx^k)\big((\pa_iN_{ik}^r)\pa_r\big)
	\]
	\[
	=-\pa_iN_{ik}^k.
	\]
	So (\ref{sum-N-ikk-is-0}) implies 
	\[
	S_i=\sum_k\hbar(\pa_i,\pa_k)=-\pa_i\sum_kN_{ik}^k=0,
	\]
	\[
	S=-\sum_{i,k=1}^n\pa_iN_{ik}^k=0.
	\]
	\qed

	\vspace{0.2in}

	\section{Type II Dual Forms  of the Squares of Nijenhuis Tensor}\label{sec-dual-form-2}
	
	In this section we express the dual forms of the squares of Nijenhuis tensor and vanishing results in forms.
	
	We consider the pair between $\Lambda^2(T_x^*M\ot\bbc)$ and
	$\Lambda^2(T_x^M\ot\bbc)$: for  $\eta\in\Lambda^2(T_x^*M\ot\bbc)$ and $u\we v\in\Lambda^2(T_x^M\ot\bbc)$,
	\[
	\om(u\we v)=\om(u,v).
	\]
	For vector $u\in T_x(M)$, let $\epsilon (u)$ be the left wedge multiplication on the exterior algebra $\Lambda(T_xM)$:
	\[
	\epsilon(u)v=u\we v\qquad \forall v\in \Lambda(T_xM\ot\bbc),
	\]
	and $\imath(u)$  the transpose of $\epsilon(u)$, that is the interior multiplication of  $\Lambda(T_x^*M\ot\bbc)$, namely,
	\[
	\om\big(\epsilon(u)v\big)=\big(\imath(u)\om\big)(v).
	\]
	
	\vspace{0.1in}
	We have the following results.
	\begin{theorem}\label{NN-thm} \ \\
		
		For (0,1)-form $\om$,
		\begin{equation}\label{01-NN-eq}
		\om\Big(N^2(X,Z; Y)\Big)=-16\big(\rhob\imath(Y)\rho\om\big)\Big( X, Z\Big)
		\end{equation}
		and
		\[
		\om\Big(JN^2(X,Z; Y)\Big)=16\imi\big(\rhob\imath(Y)\rho\om\big)\Big( X, Z\Big);
		\]

		\vspace{0.1in}
		
		For (1,0)-form $\theta$,
		\begin{equation} \label{10-NN-eq}
		\theta\Big(N^2(X, Z\we Y)\Big)=-16\big(\rho\imath(Y)\rhob\theta\big)\Big(X, Z\Big)
		\end{equation}
		and
		\[
		\theta\Big(JN^2(X,Z; Y)\Big)=-16\imi\big(\rho\imath(Y)\rhob\theta\big)\Big(X, Z\Big);
		\]
		
		\vspace{0.1in}
		
		Therefore one can extend $N^2$ and $JN^2$ to derivations on exterior algebra of differential forms.	\\
		
		We also have the following.
		
		\vspace{0.1in}
		
		For real form $\zeta$,
		\begin{equation}\label{real-NN-eq}
		\zeta\Big(N^2(X,Z; Y)\Big)
		\end{equation}
		\[
		=-8\Big\{\rhob\imath(Y)\rho(\zeta+\imi J\zeta)+\rho\imath(Y)\rhob(\zeta-\imi J\zeta)\Big\}\Big(X, Z\Big),
		\]
		and
		\[
		\zeta\Big(JN^2(X,Z; Y)\Big)
		\]
		\[
		=8\imi \Big\{\rhob\imath(Y)\rho(\zeta+\imi J\zeta)-\rho\imath(Y)\rhob(\zeta-\imi J\zeta)\Big\}\Big(X, Z\Big).
		\]

		\vspace{0.1in}
		
	\end{theorem}
	
	\proof
	For $(0,1)$ form $\om$, real vectors $X,Z$ and $Y$,  we apply Lemma \ref{N-mu-lem}. Then
	\[
	\om\Big(N^2(X,Z; Y)\Big)
	\]
	\[
	=\om\Big(N\big\{N(X, Z), Y\big\}\Big)
	\]
	\[
	=4(\rho \om)\Big(N(X, Z), Y\Big)
	\]
	\[
	=-4(\rho \om)\Big(Y, N(X, Z)\Big)
	\]
	\[
	=-4(\rho \om)\Big(\epsilon(Y) N(X, Z)\Big)
	\]
	\[
	=-4\big(\imath(Y)(\rho\om)\big)\Big( N(X, Z)\Big), 
	\]
	\[
	=-16\big(\rhob\imath(Y)\rho\om\big)\Big( N(X, Z)\Big), 
	\]
	where in the last equation we applied Lemma \ref{N-mu-lem}.
	Therefore we have
	\[
	\om\Big(N^2(X, Z; Y)\Big)=-16\big(\rhob\imath(Y)\rho\om\big)\Big( N(X, Z)\Big), 
	\]
	This equation and Lemma \ref{J-through-01-10-form-lem} imply the following equation
	\[
	\om\Big(JN^2(X, Z; Y)\Big)
	=16\imi \big(\rhob\imath(Y)\rho\om\big)\Big( N(X, Z)\Big).
	\]

	\vspace{0.1in}

	Similarly, for $(1,0))$ form $\theta$, real vectors $X,Z$ and $Y$,  we apply Lemma \ref{N-mu-lem}. Then
	\[
	\theta\Big(N^2(X, Z; Y)\Big)
	\]
	\[
	=\theta\Big(N\big\{N(X, Z), Y\big\}\Big)
	\]
	\[
	=4 (\rhob\theta)\Big(N(X, Z), Y\Big)
	\]
	\[
	=-4 (\rhob\theta)\Big(Y, N(X, Z)\Big)
	\]
	\[
	=-4 (\rhob\theta)\Big(\epsilon(Y)N(X, Z)\Big)
	\]
	\[
	=-4 \big(\imath(Y)\rhob\theta\big)\Big(N(X, Z)\Big)
	\]
	\[
	=-16\big(\rho\imath(Y)\rhob\theta\big)\Big(X, Z\Big), 
	\]
	where in the last equation we applied Lemma \ref{N-mu-lem}. Therefore we have
	\[
	\theta\Big(N^2(X, Z; Y)\Big)=-16\big(\rho\imath(Y)\rhob\theta\big)\Big(X, Z\Big).
	\]
	This equation and Lemma \ref{J-through-01-10-form-lem} imply the following equation
	\[
	\theta\Big(JN^2(X, Z; Y)\Big)=-16\imi\big(\rho\imath(Y)\rhob\theta\big)\Big(X, Z\Big).
	\]

	\vspace{0.1in}
	
	For real form $\zeta$, write $\zeta$ as the sum of (0,1) form and (1,0) form and apply
	\[
	\zeta=\frac12(\zeta+\imi J\zeta)+\frac12(\zeta-\imi J\zeta).
	\]
	and apply (\ref{01-NN-eq}) and (\ref{10-NN-eq}). Then we have
	\[
	\zeta\Big(N^2(X, Z; Y)\Big)
	\]
	\[
	=-8\Big\{\rhob\imath(Y)\rho(\zeta+\imi J\zeta)+\rho\imath(Y)\rhob(\zeta-\imi J\zeta)\Big\}\Big(X, Z\Big),
	\]
	This equation and Lemma \ref{J-through-01-10-form-lem} imply
	\[
	\zeta\Big(JN^2(X, Z; Y)\Big)
	\]
	\[
	=8\imi\Big\{\rhob\imath(Y)\rho(\zeta+\imi J\zeta)-\rho\imath(Y)\rhob(\zeta-\imi J\zeta)\Big\}\Big(X, Z\Big).
	\]

	\qed

	\vspace{0.1in}
	The following  lemma was used or will be used in several occasions in this paper.
	
	\begin{lemma}\label{J-through-01-10-form-lem}
		
		\[
		\om(JX)=-\imi\om(X)
		\]
		for	(0,1)-form $\om$;
		\[
		\theta(JX)=\imi\theta(X)
		\]
		for	(1,0)-form $\theta$.
	\end{lemma}
	
	\proof.
	We compute
	\[
	\om(JX)=\frac12\om ((JX-\imi JJX)+(JX+\imi JJX))
	\]
	\[
	=\frac12\om (JX+\imi JJX)
	\]
	\[
	=\frac12\om (JX-\imi X)
	\]
	\[
	=\frac12\om(JX)-\frac{\imi}2\om(X),
	\]
	\[
	\frac12\om(JX)=-\frac{\imi}2\om(X).
	\]
	
	Therefore
	\[
	\om(JX)=-\imi\om(X).
	\]
	
	\vspace{0.1in}
	
	Similarly, we have
	\[
	\theta(JX)=\frac12\theta((JX-\imi JJX)+(JX+\imi JJX))
	\]
	\[
	=\frac12\theta(JX-\imi JJX)
	\]
	\[
	=\frac12\theta(JX+\imi X)
	\]
	\[
	=\frac12\theta(JX)+\frac12\imi\theta(X)
	\]
	
	Therefore
	\[
	\theta(JX)=\imi\theta(X).
	\]
	
	\qed.
	
	\vspace{0.1in}

	In section \ref{sec-dual-form-1}, we defined $K,L,\hbar,\ell, S, T$ and etc. in the following we let  $X^*, Z^*, Y^*, W^*$ be the duals of real vector $X, Z, Y, W$, respectively.
	Let $\{x^i\}_{i=1}^{2m}$ be a local chart of $M$ at $x$, $\big\{\dfrac{\pa}{\pa x^i}\big\}_{i=1}^{2m}$ local frame, denote it by $\{\pa_i\}_{i=1}^{2m}$ and denote its dual
	by $\{dx^i\}_{i=1}^{2m}$. 
	
	\vspace{0.1in}
	
	We have the following dual forms.
	
	\begin{theorem}\label{L-ell-T-thm}
		We have the dual forms for $L$, $\ell$ and $T$ in \cite{Ling}.
		\[
		K(X,Z,Y,W)
		\]
		\[
		=-2\Big\{\rhob\imath(Y)\rho(W^*+\imi JW^*)+\rho\imath(Y)\rhob(W^*-\imi JW^*)\Big\}\Big(X, Z\Big),
		\]
		\[
		-2\Big\{\rhob\imath(X)\rho(W^*+\imi JW^*)+\rho\imath(Y)\rhob(W^*-\imi JW^*)\Big\}\Big(Y, Z\Big)
		\]
		\[
		-2\Big\{\rhob\imath(Y)\rho(Z^*+\imi JZ^*)+\rho\imath(Y)\rhob(Z^*-\imi JZ^*)\Big\}\Big(X, W\Big),
		\]
		\[
		-2\Big\{\rhob\imath(X)\rho(Z^*+\imi JZ^*)+\rho\imath(Y)\rhob(Z^*-\imi JZ^*)\Big\}\Big(Y, W\Big),
		\]
		
		\[
		L(X,Z,Y,W)
		\]
		\[
		=2\imi\Big\{\rhob\imath(Y)\rho(W^*+\imi JW^*)-\rho\imath(Y)\rhob(W^*-\imi JW^*)\Big\}\Big(X, Z\Big),
		\]
		\[
		+2\imi\Big\{\rhob\imath(X)\rho(W^*+\imi JW^*)-\rho\imath(Y)\rhob(W^*-\imi JW^*)\Big\}\Big(Y, Z\Big)
		\]
		\[
		+2\imi\Big\{\rhob\imath(Y)\rho(Z^*+\imi JZ^*)-\rho\imath(Y)\rhob(Z^*-\imi JZ^*)\Big\}\Big(X, W\Big),
		\]
		\[
		+2\imi\Big\{\rhob\imath(X)\rho(Z^*+\imi JZ^*)-\rho\imath(Y)\rhob(Z^*-\imi JZ^*)\Big\}\Big(Y, W\Big),
		\]
		where $X,Z,Y,W$ are real vectors and $X^*, Z^*, Y^*, W^*$ are  the duals, respectively.
		
		\vspace{0.1in}
		
		In particular,
		
		\[
		K(\pa_i,\pa_k,\pa_j,\pa_l)
		\]
		\[
		=-2\Big\{\rhob\imath(\pa_j)\rho(dx^l+\imi Jdx^l)+\rho\imath(\pa_j)\rhob(dx^l-\imi Jdx^l)\Big\}\Big(\pa_i, \pa_k\Big)
		\]
		\[
		-2\Big\{\rhob\imath(\pa_i)\rho(dx^l+\imi Jdx^l)+\rho\imath(\pa_i)\rhob(dx^l-\imi Jdx^l)\Big\}\Big(\pa_j, \pa_k\Big)
		\]
		\[
		-2\Big\{\rhob\imath(\pa_j)\rho(dx^k+\imi Jdx^k)+\rho\imath(\pa_j)\rhob(dx^k-\imi Jdx^k)\Big\}\Big(\pa_i, \pa_l\Big)
		\]
		\[
		-2\Big\{\rhob\imath(\pa_i)\rho(dx^k+\imi Jdx^k)+\rho\imath(\pa_i)\rhob(dx^k-\imi Jdx^k)\Big\}\Big(\pa_j, \pa_l\Big).
		\]
		
		\[
		L(\pa_i,\pa_k,\pa_j,\pa_l)
		\]
		\[
		=2\imi\Big\{\rhob\imath(\pa_j)\rho(dx^l+\imi Jdx^l)-\rho\imath(\pa_j)\rhob(dx^l-\imi Jdx^l)\Big\}\Big(\pa_i, \pa_k\Big)
		\]
		\[
		+2\imi\Big\{\rhob\imath(\pa_i)\rho(dx^l+\imi Jdx^l)-\rho\imath(\pa_i)\rhob(dx^l-\imi Jdx^l)\Big\}\Big(\pa_j, \pa_k\Big)
		\]
		\[
		+2\imi\Big\{\rhob\imath(\pa_j)\rho(dx^k+\imi Jdx^k)-\rho\imath(\pa_j)\rhob(dx^k-\imi Jdx^k)\Big\}\Big(\pa_i, \pa_l\Big)
		\]
		\[
		+2\imi\Big\{\rhob\imath(\pa_i)\rho(dx^k+\imi Jdx^k)-\rho\imath(\pa_i)\rhob(dx^k-\imi Jdx^k)\Big\}\Big(\pa_j, \pa_l\Big).
		\]
		
		\vspace{0.1in}
		
		\[
		\hbar(\pa_i,\pa_k)
		\]
		\[
		=-8\Big\{\rhob\imath(\pa_i)\rho(dx^k+\imi Jdx^k)+\rho\imath(\pa_i)\rhob(dx^k-\imi Jdx^k)\Big\}\Big(\pa_i, \pa_k\Big).
		\]

		\[
		\ell(\pa_i,\pa_k)
		\]
		\[
		=8\imi\Big\{\rhob\imath(\pa_i)\rho(dx^k+\imi Jdx^k)-\rho\imath(\pa_i)\rhob(dx^k-\imi Jdx^k)\Big\}\Big(\pa_i, \pa_k\Big).
		\]
		
		\vspace{0.1in}
		\[
		S_i=\sum_{k}-8\Big\{\rhob\imath(\pa_i)\rho(dx^k+\imi Jdx^k)+\rho\imath(\pa_i)\rhob(dx^k-\imi Jdx^k)\Big\}\Big(\pa_i, \pa_k\Big).
		\]
		
		\[
		S=\sum_{i,k}-8\Big\{\rhob\imath(\pa_i)\rho(dx^k+\imi Jdx^k)+\rho\imath(\pa_i)\rhob(dx^k-\imi Jdx^k)\Big\}\Big(\pa_i, \pa_k\Big).
		\]
		
		\[
		T=\sum_{i,k}8\imi\Big\{\rhob\imath(\pa_i)\rho(dx^k+\imi Jdx^k)-\rho\imath(\pa_i)\rhob(dx^k-\imi Jdx^k)\Big\}\Big(\pa_i, \pa_k\Big).
		\]
		\vspace{0.1in}
		
	\end{theorem}
	\proof
	Apply Theorem \ref{NN-thm}.
	\qed
	
	\vspace{0.1in}
	
	\begin{corollary}\label{S-T-Vanish-coro} We have the identities
		\[
		\sum_{k}\Big\{\rhob\imath(\pa_i)\rho(dx^k+\imi Jdx^k)+\rho\imath(\pa_i)\rhob(dx^k-\imi Jdx^k)\Big\}\Big(\pa_i, \pa_k\Big)\equiv 0
		\]
		for $1\le i\le n$,
		\[
		\sum_{i,k}\Big\{\rhob\imath(\pa_i)\rho(dx^k+\imi Jdx^k)\Big\}\Big(\pa_i, \pa_k\Big)\equiv 0,
		\]
		\[
		\sum_{i,k}\Big\{\rho\imath(\pa_i)\rhob(dx^k-\imi Jdx^k)\Big\}\Big(\pa_i, \pa_k\Big)\equiv 0.
		\]	
		
	\end{corollary}
	\proof The dual forms of $S_i$, $S$ and  $T$ in the above theorem and
	Theorem \ref{S-vanishes-thm} and Theorem 2.1 in \cite{Ling} imply the following.
	\[
	\sum_{k}\Big\{\rhob\imath(\pa_i)\rho(dx^k+\imi Jdx^k)+\rho\imath(\pa_i)\rhob(dx^k-\imi Jdx^k)\Big\}\Big(\pa_i, \pa_k\Big)\equiv 0
	\]
	
	\[
	\sum_{i,k}\Big\{\rhob\imath(\pa_i)\rho(dx^k+\imi Jdx^k)+\rho\imath(\pa_i)\rhob(dx^k-\imi Jdx^k)\Big\}\Big(\pa_i, \pa_k\Big)\equiv 0,
	\]	
	and
	\[
	\sum_{i,k}\Big\{\rhob\imath(\pa_i)\rho(dx^k+\imi Jdx^k)-\rho\imath(\pa_i)\rhob(dx^k-\imi Jdx^k)\Big\}\Big(\pa_i, \pa_k\Big)
	\equiv 0.
	\]
	From the last two equations we have
	\[
	\sum_{i,k}\Big\{\rhob\imath(\pa_i)\rho(dx^k+\imi Jdx^k)\Big\}\Big(\pa_i, \pa_k\Big)\equiv 0,
	\]
	\[
	\sum_{i,k}\Big\{\rho\imath(\pa_i)\rhob(dx^k-\imi Jdx^k)\Big\}\Big(\pa_i, \pa_k\Big)\equiv 0.
	\]		
	\qed
	
	\vspace{0.1in}
	
	We have results from  \cite{Ling}, for example, we have the following: 
	$\ell\equiv0$ implies $N^2\equiv 0$.
	
	\vspace{0.1in}
	
	Remark. 
	The formulations of $N^2$, $L$, $\ell$, $T$, $T=0$ above actually show that
	they are metric independent those they are defined with some metric in \cite{Ling}.
	
	\vspace{0.2in}

	\begin{center}
		Acknowledgement
	\end{center}
	The author is grateful to H. Lin, Z. Lu, L. Ni and S. Wilson for communications.

\end{document}